\documentclass[11pt]{amsart}
\usepackage[T1]{fontenc}
\usepackage{lmodern}
\usepackage{amsmath,amssymb,mathtools,amsfonts}
\usepackage{amsthm}
\usepackage{booktabs}
\usepackage{mathrsfs}
\usepackage{xcolor}
\usepackage{microtype}
\usepackage{hyperref}
\usepackage{url}
\usepackage{enumitem}
\hypersetup{colorlinks=true,linkcolor=blue!50!black,citecolor=blue!50!black,urlcolor=blue!50!black,pdftitle={The covering radius of the first-order Reed--Muller code in nine variables is 242}}
\numberwithin{equation}{section}
\allowdisplaybreaks[2]
\theoremstyle{plain}
\newtheorem{theorem}{Theorem}[section]
\newtheorem{proposition}[theorem]{Proposition}
\newtheorem{lemma}[theorem]{Lemma}

\theoremstyle{definition}

\newtheorem{certificate}[theorem]{Finite certificate}
\theoremstyle{remark}

\newcommand{\F}{\mathbb F}
\newcommand{\V}{\F_2^9}
\newcommand{\NL}{\operatorname{NL}}
\newcommand{\wt}{\operatorname{wt}}
\newcommand{\RM}{\operatorname{RM}}
\newcommand{\supp}{\operatorname{supp}}
\newcommand{\one}{\mathbf{1}}
\newcommand{\modtwo}{\pmod 2}
\title[The covering radius of $\RM(1,9)$ is 242]{The covering radius of the first-order Reed--Muller code in nine variables is 242}
\author{Yi Li}
\address{School of Computer Science, Shanghai Jiao Tong University, Shanghai, China}
\date{October 31, 2026}

\begin{document}
\begin{abstract}
We prove that the covering radius of the binary first-order Reed--Muller
code $\RM(1,9)$ is exactly $242$.  Equivalently, the maximum nonlinearity
of a Boolean function in nine variables is $242$.  Our proof combines a
one-point rigidity theorem with a spectral and computational exclusion of
the two larger candidate nonlinearities.  We first show that every
nine-variable Boolean function of nonlinearity $243$ admits a unique
one-point modification that raises its nonlinearity to $244$: the point
to be changed is the vector sum of its support.  This part uses exact
Walsh-update identities, integer interpolation, rigidity of a quadratic
residue, and a sharp Walsh $\ell^1$ estimate for cubic Boolean functions.
To exclude nonlinearity $244$, we use Parseval's identity and congruences
of integral Walsh coefficients to force an exceptionally restricted
spectrum.  An associated cubic residue phase must attain the extremal
Walsh $\ell^1$ bound $10240$.  We eliminate one equality spectrum by
hyperplane slicing; the other forces the original Boolean function from
the signs of the cubic phase and yields a double-complement-pairing
identity for its cubic top part.  The classification of alternating
trilinear forms over $\F_2$, an exact filter of their representatives,
and two exhaustive enumerations of $2^{18}$ quadratic supplements then
rule out every remaining candidate.  Together with the known
construction attaining nonlinearity $242$, these arguments establish
the asserted covering radius.
\end{abstract}
\maketitle

\section{Introduction}
A \emph{Boolean function} in $n$ variables is a map
$f:\F_2^n\to\F_2$.  We identify $f$ with its truth table, a binary word
of length $2^n$.  For two Boolean functions $f$ and $g$, their Hamming
distance $d_H(f,g)$ is the number of inputs on which their values differ.
The \emph{nonlinearity} of $f$ is the minimum distance from $f$ to an
affine Boolean function:
\begin{equation}\label{eq:intro-nonlinearity}
 \NL(f)
 =\min_{\substack{a\in\F_2^n\\ b\in\F_2}}
       d_H\bigl(f,a\cdot x+b\bigr)
 =\min_{\substack{a\in\F_2^n\\ b\in\F_2}}
       \wt\bigl(f(x)+a\cdot x+b\bigr).
\end{equation}
Here $\wt$ denotes Hamming weight, and addition inside the weight is
in $\F_2$.  Thus nonlinearity measures how far a Boolean function lies
from every affine approximation.  Large nonlinearity is a central
criterion in the design of cryptographic Boolean functions because it
limits their correlation with affine functions; it is equally natural
from the standpoint of coding theory \cite{Carlet2021}.

For a binary code $\mathcal C\subseteq\F_2^N$, its \emph{covering
radius} is the largest distance of a binary word from its nearest
codeword:
\begin{equation}\label{eq:intro-covering-definition}
 \rho(\mathcal C)
   =\max_{u\in\F_2^N}\ \min_{c\in\mathcal C}d_H(u,c).
\end{equation}
Equivalently, $\rho(\mathcal C)$ is the smallest integer $R$ such that
the Hamming balls of radius $R$ centred at the codewords cover the entire
ambient space $\F_2^N$.  The binary first-order Reed--Muller code
\[
 \RM(1,n)
  =\bigl\{\bigl(a\cdot x+b\bigr)_{x\in\F_2^n}:
          a\in\F_2^n,\ b\in\F_2\bigr\}
  \subseteq\F_2^{2^n}
\]
consists precisely of the truth tables of affine Boolean functions.
Every word in $\F_2^{2^n}$ is the truth table of exactly one Boolean
function on $\F_2^n$.  Consequently, substituting this description of
$\RM(1,n)$ into the covering-radius definition and comparing with
\eqref{eq:intro-nonlinearity} gives the fundamental identity
\begin{equation}\label{eq:intro-radius}
 \boxed{\rho\bigl(\RM(1,n)\bigr)
   =\max_{f:\F_2^n\to\F_2}\NL(f).}
\end{equation}
In other words, determining the covering radius of $\RM(1,n)$ and
determining the maximum possible nonlinearity of an $n$-variable
Boolean function are \emph{exactly the same extremal problem}.

The difficulty of this problem is particularly pronounced for odd
$n$.  The Walsh transform
\[
 W_f(a)=\sum_{x\in\F_2^n}(-1)^{f(x)+a\cdot x}
\]
provides the equivalent spectral expression
\begin{equation}\label{eq:intro-walsh-expression}
 \NL(f)=2^{n-1}-\frac12\max_{a\in\F_2^n}|W_f(a)|.
\end{equation}
When $n$ is even, bent functions achieve a perfectly flat spectrum,
providing an extremal benchmark.  No such spectrum is possible when
$n$ is odd, since Parseval's identity would require all Walsh
magnitudes to be the nonintegral number $2^{n/2}$.  The resulting
optimization is global: a large nonlinearity requires the entire Walsh
spectrum, rather than one or a few coefficients, to remain small.

The case $n=9$ illustrates how hard this global constraint can be.
Appending one variable to an eight-variable bent function yields
nonlinearity $240$; constructions of Kavut, Maitra and Y\"ucel
\cite{KavutMaitraYucel2007} attained $241$, and Kavut and Y\"ucel
\cite{KavutYucel2010} obtained functions with nonlinearity $242$.
On the other hand, Walsh integrality and Parseval's identity imply the
upper bound $244$.  The long-standing uncertainty was therefore
confined to the exceptionally narrow interval
\begin{equation}\label{eq:intro-bounds}
 242\le \rho\bigl(\RM(1,9)\bigr)\le244.
\end{equation}
Closing this gap is not a small-scale search problem.  There are
$2^{512}$ nine-variable truth tables; changing one output bit affects
all $512$ Walsh coefficients at once.

Our principal result settles the problem.
\begin{theorem}[Main theorem]\label{thm:main}
The covering radius of the binary first-order Reed--Muller code
$\RM(1,9)$ is
\[
 \boxed{\rho\bigl(\RM(1,9)\bigr)=242.}
\]
Equivalently, every nine-variable Boolean function has nonlinearity
at most $242$, and this value is attained.
\end{theorem}

A key feature of our approach is a precise connection between the two
potential improvements beyond $242$.  For $t\in\V$, let $\delta_t$ be
the indicator of the singleton $\{t\}$, and define the \emph{support
sum} of a Boolean function by
\[
 p(f)=\sum_{x\in\supp(f)}x\in\V,
\]
where the sum is taken in the vector space $\V=\F_2^9$.
\begin{theorem}[One-point rigidity]\label{thm:flip}
If $\NL(f)=243$, then
\[
 \NL\bigl(f+\delta_{p(f)}\bigr)=244,
\]
and $p(f)$ is the unique point whose flip increases nonlinearity to
$244$.  Conversely, if $\NL(g)=244$, then
$\NL(g+\delta_t)=243$ for every $t\in\V$ and $p(g)=0$.
Consequently the map
\[
 \{g:\NL(g)=244\}\times\V\longrightarrow\{f:\NL(f)=243\},
 \qquad(g,t)\longmapsto g+\delta_t
\]
is a bijection.
\end{theorem}

We prove the one-point theorem by exploiting the integer structure of
Walsh coefficients.  Exact interpolation shows that successive binary
residues of a normalized Walsh spectrum are Boolean polynomials of
small degree.  A quadratic residue rigidity argument and a sharp
$\ell^1$ estimate for cubic spectra then identify the unique improving
flip.  It follows that nonlinearity $243$ cannot occur unless
nonlinearity $244$ does; hence it remains to exclude the latter.

Suppose, therefore, that a Boolean function of nonlinearity $244$
exists.  Parseval's identity and integral Walsh congruences force
Walsh magnitudes $8$ and $24$, occurring $64$ and $448$ times,
respectively.  The next binary residue yields a cubic Boolean phase
$h$ for which $\|W_h\|_1\ge10240$.  Our cubic spectral inequality
shows that equality must hold and leaves only two possible absolute
spectra.  We eliminate one by comparing the Walsh energies on the two
sides of an eight-dimensional hyperplane.  In the remaining spectrum,
the original function is determined by the signs of $W_h$; comparison
of its degree-six coefficients with the cubic top part $\omega$ of $h$
gives a nonlinear double-complement-pairing identity
$\omega^{\#\#}=\omega$.  The classification of nine-dimensional
alternating trilinear forms \cite{BrierLangevin2003,HoraPudlak2021}
combined with an exact filter of the fixed-point identity, then
leaves only two cubic representatives.  For each of them, nine linear face constraints and nine phase-translation directions
reduce the $36$ quadratic coefficients to a complete $18$-dimensional
transversal.  Two exact $2^{18}$-state enumerations exclude every
remaining equality-spectrum candidate.

We obtained the conclusion
$\rho\bigl(\RM(1,9)\bigr)=242$ on September 30, 2026.
The present paper gives our independent proof of this result.
Contemporaneous results were also obtained by Wang, Lou and Wang
\cite{WangLouWang2026} and by Hou, Sabo and Slaughter
\cite{HouSaboSlaughter2026}.  Our approach combines one-point
rigidity, arithmetic restrictions on integral Walsh spectra,
and sharp estimates for cubic Boolean functions.  The classification
of alternating trilinear forms reduces the remaining obstruction
to two cubic representatives; linear constraints and translation
symmetries then leave $2^{18}$ quadratic supplements for each.
An exhaustive computer calculation eliminates these finite
possibilities and completes the proof.

\section{Walsh and integer-polynomial preliminaries}
Set $V=\F_2^9$ and $[9]=\{1,\ldots,9\}$.  For $S\subseteq[9]$ write
$x^S=\prod_{i\in S}x_i$, with $x^\varnothing=1$, and write $S^c=[9]\setminus S$.
The algebraic normal form (ANF) of a Boolean function is its unique
squarefree polynomial over $\F_2$.  Integer sums below are explicitly
distinguished from sums in $\F_2$.

For a Boolean function $f$, define its sign and unnormalized Walsh
transform by
\[
 s_f(x)=(-1)^{f(x)},\qquad
 W_f(a)=\sum_{x\in V}(-1)^{f(x)+a\cdot x},\qquad a\in V.
\]
Character orthogonality gives
\begin{align}
 \sum_{a\in V}(-1)^{a\cdot x}W_f(a)&=512s_f(x),
       \label{eq:inversion}\\
 \sum_{a\in V}W_f(a)^2&=512^2=262144,\label{eq:parseval}\\
 \NL(f)&=256-\tfrac12\|W_f\|_\infty.\label{eq:nl-w}
\end{align}
In particular $\|W_f\|_\infty\ge 24$ because this maximum is an even
integer and $\sqrt{512}>22$.  This proves the upper bound $244$
in~\eqref{eq:intro-bounds}.  Symmetry of the Walsh--Hadamard matrix also yields
\begin{equation}\label{eq:dual-pairing}
 \sum_{a\in V}s_h(a)W_f(a)=\sum_{x\in V}s_f(x)W_h(x).
\end{equation}
For any $t\in V$, a direct change in one sign gives the exact update
\begin{equation}\label{eq:flip-formula}
 W_{f+\delta_t}(a)=W_f(a)-2s_f(t)(-1)^{a\cdot t}.
\end{equation}
It follows that $|\NL(f+\delta_t)-\NL(f)|\le1$.

\subsection{An integer Walsh polynomial}
Every integer-valued function on $\{0,1\}^9$ has a unique integer
multilinear representation, by M\"obius inversion:
\begin{equation}\label{eq:mobius}
 T(a)=\sum_{S\subseteq[9]}T_Sa^S,\qquad
 T_S=\sum_{U\subseteq S}(-1)^{|S|-|U|}T(\one_U).
\end{equation}
In particular, if $m$ divides every value of $T$, then it divides
every coefficient $T_S$.  This fact justifies division in $\mathbb Z$
\emph{before} reduction modulo two.

For $S\subseteq[9]$, let
\[
 M_S(f)=\sum_{\substack{x\in V\\x_i=1\;(i\in S)}} f(x)
\]
be the integer weight on the indicated coordinate face.  Expanding
$(-1)^{a\cdot x}=\prod_{i:x_i=1}(1-2a_i)$ shows that the coefficient
$C_S(f)$ of $a^S$ in $W_f(a)$ is
\begin{equation}\label{eq:face-walsh}
 C_S(f)=(-1)^{|S|}512-2(-2)^{|S|}M_S(f).
\end{equation}
For $S=\varnothing$ this reads $W_f(0)=512-2\wt(f)$.
For $|S|\ge 3$ the right-hand side is divisible by $16$, and for
$|S|\ge4$ it is divisible by $32$.

\begin{lemma}[Integral residue degrees]\label{lem:residue}
If $W_f(a)\equiv c\pmod{8}$ at every frequency, then
$((W_f-c)/8)\bmod2$ has degree at most two.  If
$W_f(a)\equiv c\pmod{16}$, then $((W_f-c)/16)\bmod2$ has degree
at most three.
\end{lemma}
\begin{proof}
Both quotients are integer-valued, so their coefficients are obtained
by exact division of the integer coefficients in~\eqref{eq:mobius}.
Equation~\eqref{eq:face-walsh} makes all coefficients of degree at least
three even after division by $8$, and all coefficients of degree at least
four even after division by $16$.  Reducing the already integral
multilinear representation modulo two proves the assertions.
\end{proof}

We recall the standard Reed--Muller minimum-weight estimate and give a proof for completeness.
\begin{lemma}[Reed--Muller minimum weight]\label{lem:minweight}
Let $m\ge1$ and $0\le d\le m$.  If
$u:\F_2^m\to\F_2$ is nonzero and $\deg u\le d$, then
$\wt(u)\ge2^{m-d}$.
\end{lemma}
\begin{proof}
We argue by induction on the number $m\ge1$ of variables.
For $m=1$, a nonzero Boolean function is either the constant
function $1$, of weight $2$, or a nonconstant affine function, of
weight $1$.  In both cases, the asserted inequality holds for
all allowed values of $d$.

Now let $m\ge2$ and assume the result for $m-1$ variables.
If $d=0$, then $u$ is the nonzero constant function and
$\wt(u)=2^m=2^{m-d}$.  We may therefore take $d\ge1$.
Write $x=(x',z)$, with $x'\in\F_2^{m-1}$ and $z\in\F_2$.
Separating the terms of the algebraic normal form according to
whether they contain $z$ gives
\[
 u(x',z)=u_0(x')+z\,u_1(x'),
 \qquad \deg u_0\le d,\quad \deg u_1\le d-1,
\]
where addition is in $\F_2$.  The restrictions to $z=0$ and $z=1$
are $u_0$ and $u_0+u_1$, respectively.  In particular,
\[
 \wt(u)=\wt(u_0)+\wt(u_0+u_1).
\]

If $u_1=0$, the two restrictions coincide.  Since $u$ is nonzero,
so is $u_0$, and hence $\wt(u)=2\wt(u_0)$.  Apply the induction
hypothesis to $u_0$ with degree bound
$d'=\min\{d,m-1\}$.  It follows that
\[
 \wt(u)\ge 2\cdot2^{(m-1)-d'}
            =2^{m-d'}\ge2^{m-d}.
\]

Now suppose $u_1\ne0$.  For every $x'$ with $u_1(x')=1$, the two
values $u(x',0)=u_0(x')$ and $u(x',1)=u_0(x')+1$ are different;
therefore exactly one of them equals $1$.  Such an $x'$ contributes
one to $\wt(u)$, irrespective of $u_0(x')$.  Consequently,
\[
 \wt(u)\ge\wt(u_1).
\]
Since $u_1$ is a nonzero function in $m-1$ variables of degree at
most $d-1$, the induction hypothesis yields
\[
 \wt(u)\ge\wt(u_1)
       \ge2^{(m-1)-(d-1)}=2^{m-d}.
\]
Both cases establish the claim.
\end{proof}

Recall that if a Boolean function $q$ has degree at most two, its
quadratic polar form has even rank $r\le8$, and the nonzero entries of
$W_q$ have magnitude $2^{9-r/2}$ at precisely $2^r$ frequencies.
Consequently
\begin{equation}\label{eq:quadratic-l1}
 \|W_q\|_1=2^{9+r/2}\le8192.
\end{equation}
This follows by squaring $W_q$ and summing over the radical of the
polar form; alternatively it is the standard quadratic Walsh formula
\cite{Carlet2021}.

\section{The sharp Walsh norm bound for cubics}
The following lemma supplies the common analytic input for both
remaining nonlinearities.

\begin{lemma}[Exact weight expansion]\label{lem:weight-expansion}
Let $u:\F_2^m\to\F_2$ be the XOR of its distinct ANF monomials with
supports in a set $\mathscr M$, including the empty support when
appropriate.  Then, as an integer identity,
\begin{equation}\label{eq:weight-expansion}
 \wt(u)=\sum_{\varnothing\ne\mathscr A\subseteq\mathscr M}
             (-2)^{|\mathscr A|-1}2^{m-|\bigcup_{J\in\mathscr A}J|}.
\end{equation}
If $m\ge7$ and $\deg u\le3$, then $4\mid\wt(u)$ and
$8\mid W_u(a)$ at every frequency.  If $m\ge4$ and $\deg u\le3$,
then $\wt(u)$ is even.
\end{lemma}
\begin{proof}
Write the algebraic normal form of $u$ as
\[
 u(x)=\sum_{J\in\mathscr M}x^J\pmod2,
 \qquad x^J=\prod_{j\in J}x_j.
\]
For a Boolean input $x$, each monomial $x^J$ is either $0$ or $1$.
Thus $(-1)^{x^J}=1-2x^J$, and taking the sign of the XOR of the
monomials gives the following identity of \emph{ordinary integers}:
\[
 (-1)^{u(x)}
  =\prod_{J\in\mathscr M}(-1)^{x^J}
  =\prod_{J\in\mathscr M}(1-2x^J).
\]
Expanding the right-hand product amounts to choosing, for each
$J\in\mathscr M$, either $1$ or $-2x^J$.  Consequently,
\[
 (-1)^{u(x)}
 =1+\sum_{\varnothing\ne\mathscr A\subseteq\mathscr M}
       (-2)^{|\mathscr A|}\prod_{J\in\mathscr A}x^J.
\]
On $\{0,1\}^m$ we have $x_i^2=x_i$, so
$\prod_{J\in\mathscr A}x^J=x^{\bigcup_{J\in\mathscr A}J}$.
Moreover, the Boolean value $u(x)$, viewed as the integer $0$ or
$1$, satisfies $u(x)=(1-(-1)^{u(x)})/2$.  Substituting the preceding
expansion therefore yields the exact integer identity
\[
 u(x)=\sum_{\varnothing\ne\mathscr A\subseteq\mathscr M}
       (-2)^{|\mathscr A|-1}
       x^{\bigcup_{J\in\mathscr A}J}.
\]
For any $U\subseteq[m]$, precisely $2^{m-|U|}$ inputs have
$x_i=1$ for all $i\in U$; hence
$\sum_{x\in\F_2^m}x^U=2^{m-|U|}$.  Summing the last identity
over all $x$ now gives~\eqref{eq:weight-expansion}.

Suppose next that $m\ge7$ and $\deg u\le3$, so $|J|\le3$ for
all $J\in\mathscr M$.  If $|\mathscr A|=1$, its contribution to
\eqref{eq:weight-expansion} is divisible by $2^{m-3}$; if
$|\mathscr A|=2$, the union has at most six elements and the
coefficient $-2$ supplies an additional factor of two, giving
divisibility by $2^{m-5}$.  Both quantities are multiples of four
when $m\ge7$.  For $|\mathscr A|\ge3$, the coefficient
$(-2)^{|\mathscr A|-1}$ itself is divisible by four.  Hence
$4\mid\wt(u)$.  Since $u(x)+a\cdot x$ also has degree at most
three, the same divisibility applies to its weight, and
\[
 W_u(a)=2^m-2\wt(u+a\cdot x)\equiv0\pmod8.
\]
Finally, for $m\ge4$ the parity of the weight is
$\sum_x u(x)\pmod2$.  In the ANF, a monomial $x^J$ contributes
$2^{m-|J|}$ to this sum.  As $|J|\le3<m$, every contribution is
even.  Thus $\wt(u)$ is even.
\end{proof}

\begin{lemma}[A shared residue and an affine next digit]\label{lem:affine-digit}
Let $h:\V\to\F_2$ have degree at most three and put $Q=W_h$.
All $Q(a)$ have a common residue modulo $16$, which is either $0$
or $8$.  Moreover, for distinct $i,j$,
\begin{equation}\label{eq:mixed-differences}
 16\mid\bigl(Q(a+e_i)-Q(a)\bigr),\qquad
 32\mid\Delta_i\Delta_jQ(a).
\end{equation}
Consequently $Q/16\bmod2$ is affine when $16\mid Q$, and
$(Q-8)/16\bmod2$ is affine when $Q\equiv8\pmod{16}$.
\end{lemma}
\begin{proof}
For an integer-valued function $T$ on $V$, define
$\Delta_iT(a)=T(a+e_i)-T(a)$, where the subtraction is performed in
$\mathbb Z$.  By comparing the terms in $Q(a+e_i)$ and $Q(a)$,
\[
 \Delta_iQ(a)
  =-2\sum_{x_i=1}(-1)^{h(x)+a\cdot x}.
\]
Applying a second difference in a distinct direction $e_j$ gives
\[
 \begin{split}
 \Delta_i\Delta_jQ(a)
 &=Q(a+e_i+e_j)-Q(a+e_i)-Q(a+e_j)+Q(a)\\
 &=4\sum_{x_i=x_j=1}(-1)^{h(x)+a\cdot x}.
 \end{split}
\]
Fixing one coordinate of $h$ leaves a Boolean polynomial of degree
at most three in eight variables; fixing two distinct coordinates
leaves one in seven variables.  The signs contributed by the fixed
coordinates are immaterial to divisibility.  Lemma~\ref{lem:weight-expansion}
shows that the two sums on the right are divisible by eight.
The factors $-2$ and $4$ therefore give, respectively,
$16\mid\Delta_iQ(a)$ and $32\mid\Delta_i\Delta_jQ(a)$.

Any two points of $V$ can be joined by changes of one coordinate
at a time.  Since each such change alters $Q$ by a multiple of
$16$, the value of $Q(a)$ is independent of $a$ modulo $16$.
Lemma~\ref{lem:weight-expansion}, applied in nine variables, also
shows that every $Q(a)$ is divisible by $8$.  Its common residue
modulo $16$ is therefore some $r\in\{0,8\}$.

Define the integer-valued function and its binary reduction by
\[
 T(a)=\frac{Q(a)-r}{16}\in\mathbb Z,
 \qquad b(a)=T(a)\bmod2\in\F_2.
\]
The constant $r$ disappears under a mixed difference.  For distinct
$i,j$, the divisibility proved above implies
\[
 \Delta_i\Delta_jT(a)
    =\frac{\Delta_i\Delta_jQ(a)}{16}\equiv0\pmod2.
\]
Reduction modulo two turns subtraction into addition.  Therefore
$b$ satisfies
\[
 \begin{split}
 D_iD_jb(a)
 &=b(a+e_i+e_j)+b(a+e_i)\\
 &\quad+b(a+e_j)+b(a)=0,
 \end{split}
 \qquad i\ne j,
\]
where $D_i$ denotes the usual Boolean derivative in direction $e_i$.
To see what this implies, write the ANF of $b$ as
$b(a)=\sum_{S\subseteq[9]}b_Sa^S$.  For every squarefree monomial,
\[
 D_iD_j(a^S)=
 \begin{cases}
 a^{S\setminus\{i,j\}},&\{i,j\}\subseteq S,\\
 0,&\text{otherwise}.
 \end{cases}
\]
For fixed $i,j$, distinct supports containing $\{i,j\}$ give
distinct monomials after removing those two indices, so their
coefficients cannot cancel.  If $b_S=1$ for any $|S|\ge2$, choose
distinct $i,j\in S$; the coefficient of $a^{S\setminus\{i,j\}}$
in $D_iD_jb$ would be $1$, contrary to the displayed identity.
Thus $b_S=0$ whenever $|S|\ge2$, and
$b(a)=b_\varnothing+\sum_{i=1}^9b_{\{i\}}a_i$ is affine.
Taking $r=0$ or $r=8$ gives exactly the two conclusions of the
lemma.
\end{proof}

\begin{theorem}[Cubic spectral rigidity]\label{thm:cubic}
For every nine-variable Boolean function $h$ of degree at most three,
\begin{equation}\label{eq:cubic-l1}
 \boxed{\|W_h\|_1\le10240.}
\end{equation}
If $W_h\equiv8\pmod{16}$, then the sharper bound
$\|W_h\|_1\le10048$ holds.  Equality in~\eqref{eq:cubic-l1}
is possible only with one of the following absolute spectra:
\begin{alignat}{2}
 \mathrm{(A)}\quad& \#\{|W_h|=16\}=448,\quad
                 &&\#\{|W_h|=48\}=64;\label{eq:spectrum-A}\\
 \mathrm{(B)}\quad& \#\{|W_h|=0\}=64,\quad
                 &&\#\{|W_h|=16\}=256,\quad\#\{|W_h|=32\}=192.
                 \label{eq:spectrum-B}
\end{alignat}
\end{theorem}
\begin{proof}
Write $Q=W_h$.  By Lemma~\ref{lem:affine-digit}, either $16\mid Q$
or $Q\equiv8\pmod{16}$.

First assume $16\mid Q$ and put $v=Q/16$.  The parity $r=v\bmod2$
is affine, so it is either constantly zero, constantly one, or balanced.
Parseval gives
\begin{equation}\label{eq:cubic-energy}
 \sum_a v(a)^2=1024.
\end{equation}
For even integers $z$ and odd integers $w$, respectively, the exact
inequalities
\begin{equation}\label{eq:integer-inequalities}
 |z|\le\frac{z^2}{2},\qquad
 |w|\le\frac{w^2+3}{4},\qquad
 |w|\le\frac{w^2+1}{2}
\end{equation}
follow by factoring the nonnegative differences.  If $r=0$, the first
bound yields $\|Q\|_1\le8192$.  If $r=1$, the second yields
\[
 \|Q\|_1\le16\frac{1024+3\cdot512}{4}=10240.
\]
Equality requires $|v|\in\{1,3\}$ pointwise.  The energy equation
then fixes $448$ entries of magnitude one and $64$ of magnitude three,
which is spectrum~\eqref{eq:spectrum-A}.

If $r$ is balanced, there are $256$ odd positions.  Apply the last
inequality in~\eqref{eq:integer-inequalities} there and the first one
on the remaining positions.  We obtain
\[
 \sum_a|v(a)|\le\frac{1024+256}{2}=640.
\]
Equality forces the odd entries to be $\pm1$ and the even entries
$0$ or $\pm2$.  The energy of the odd half is $256$, leaving
$768=192\cdot4$, so there are $192$ even entries of magnitude two
and $64$ zeros.  This is spectrum~\eqref{eq:spectrum-B}.

Finally suppose $Q\equiv8\pmod{16}$.  The odd-valued function
$u=Q/8$ has $\lambda=(u-1)/2\bmod2$ affine.  Define
$z(a)=(-1)^{\lambda(a)}u(a)$.  Then $z\equiv1\pmod4$ and
$|z|=|u|$.  Because $\lambda=c+b\cdot a$, Walsh inversion and
Parseval give
\begin{equation}\label{eq:z-moments}
 \sum_a z(a)=64(-1)^c s_h(b)\in\{-64,64\},\qquad
 \sum_a z(a)^2=4096.
\end{equation}
For every integer $z\equiv1\pmod4$,
\begin{equation}\label{eq:majorant}
 16|z|\le3z^2-2z+15.
\end{equation}
For $z\ge1$ the difference of right and left sides is
$3(z-1)(z-5)$; for $z\le-3$ it equals $(3z+5)(z+3)$.
Both are nonnegative on the allowed progression.  Summing and using
\eqref{eq:z-moments} gives
\[
 \|Q\|_1=8\sum_a|z(a)|
 \le\frac{3\cdot4096+128+15\cdot512}{2}=10048.
\]
This is strictly smaller than $10240$ and finishes the proof.
\end{proof}

\section{Excluding the nonlinearity-243 case}
\subsection{Parity and the support sum}
Every affine Boolean function on $V$ has even weight.  Thus all
Hamming distances from $f$ to affine functions have the parity of
$\wt(f)$, and
\begin{equation}\label{eq:nl-parity}
 \NL(f)\equiv\wt(f)\modtwo.
\end{equation}
The parity of $\wt(f)$ is also the coefficient of $x_1\cdots x_9$
in the ANF.  Consequently $\NL(f)=243$ forces
\begin{equation}\label{eq:odd-degree}
 \wt(f)\equiv1\modtwo,\qquad\deg(f)=9.
\end{equation}
If $S=\supp(f)$ and $k_a=\#\{x\in S:a\cdot x=1\}$,
then $k_a\equiv a\cdot p(f)\pmod2$.  Expanding the character sum
therefore gives the useful congruence
\begin{equation}\label{eq:support-sum}
 W_f(a)\equiv -2\wt(f)+4(a\cdot p(f))\pmod8,
\end{equation}
where $a\cdot p(f)$ is represented as $0$ or $1$ on the right.

Suppose $\NL(f)=243$.  Translating the input by $p(f)$ makes the
support sum zero: since $\wt(f)$ is odd,
\[
 p(f(\,\cdot\,+p(f)))=p(f)+\wt(f)p(f)=0.
\]
The common Walsh residue modulo eight is then $2$ or $6$; output
complementation, if necessary, changes it to $2$.  These operations
preserve nonlinearity and the absolute Walsh spectrum.  We may
therefore assume throughout the next two subsections that
\begin{equation}\label{eq:normalized-243}
 p(f)=0,\qquad W_f(a)\equiv2\pmod8,
 \qquad\|W_f\|_\infty=26.
\end{equation}
No assumption about the unknown sign $s_f(0)$ is made.

\subsection{A rigid quadratic residue}
From~\eqref{eq:normalized-243},
\[
 W_f(a)\in\{-22,-14,-6,2,10,18,26\}.
\]
Put $v=(W_f-2)/8$, so that
$v(a)\in\{-3,-2,-1,0,1,2,3\}$, and set
$q=v\bmod2$.  Lemma~\ref{lem:residue} gives $\deg q\le2$.
Writing $\sigma=s_f(0)\in\{1,-1\}$, inversion at zero and Parseval
produce the exact moment equations
\begin{equation}\label{eq:243-moments}
 \sum_a v(a)=64\sigma-128,\qquad
 \sum_a v(a)^2=4128-32\sigma\ge4096.
\end{equation}
Let $m=\wt(q)$.  On a zero of $q$, $v$ is even and $v^2\le4$;
on a one, $v$ is odd and $v^2\le9$.  Hence
\begin{equation}\label{eq:243-lower-m}
 4096\le\sum_a v(a)^2\le4(512-m)+9m=2048+5m,
 \qquad m\ge410.
\end{equation}
If $q$ were not identically one, $1+q$ would be a nonzero Boolean
function of degree at most two.  Lemma~\ref{lem:minweight} would give
$\wt(1+q)\ge128$, or $m\le384$, a contradiction.
We have proved
\begin{equation}\label{eq:243-four}
 q\equiv1,\quad W_f(a)\equiv10\pmod{16},\quad
 W_f(a)\in\{-22,-6,10,26\}.
\end{equation}

\subsection{Determining the sign at the origin}
The decisive point is that the normalization actually forces
$s_f(0)=1$.  Suppose, to the contrary, that $\sigma=-1$.
Let $A,B,C,D$ be the respective multiplicities of the four values
$-22,-6,10,26$ in~\eqref{eq:243-four}.  Counting, inversion at the
origin, and Parseval give
\begin{align*}
 A+B+C+D&=512,\\
 -22A-6B+10C+26D&=-512,\\
 484A+36B+100C+676D&=262144.
\end{align*}
Elementary elimination yields
\begin{equation}\label{eq:ABCD}
 A+D=456,\qquad B=808-3A,\quad C=3A-752,\quad
 251\le A\le269.
\end{equation}
Define a Boolean phase on frequency space by the \emph{integer}
quotient
\begin{equation}\label{eq:phase-243}
 h(a)=\frac{W_f(a)-10}{16}\pmod2.
\end{equation}
Lemma~\ref{lem:residue} gives $\deg h\le3$.
Here the \emph{sign of the phase} at a frequency $a$ means
$s_h(a)=(-1)^{h(a)}\in\{1,-1\}$; it is not, in general, the
ordinary sign of the integer $W_f(a)$.  To determine $s_h(a)$,
first evaluate the integer quotient in~\eqref{eq:phase-243} and
then reduce it modulo two.  For the four possible Walsh values,
this gives
\[
\begin{array}{c|rrrr}
 W_f(a)&-22&-6&10&26\\ \hline
 (W_f(a)-10)/16&-2&-1&0&1\\
 h(a)&0&1&0&1\\
 s_h(a)=(-1)^{h(a)}&1&-1&1&-1
\end{array}
\]
Indeed, the quotients $-2$ and $0$ are even, so the corresponding
phase values are zero and their signs are $1$.  The quotients
$-1$ and $1$ are odd, so the phase values are one and their signs
are $-1$.  This also illustrates why reduction modulo two must be
performed \emph{after} the integer division.  Thus
\begin{align}
 \wt(h)&=B+D=1264-4A,\label{eq:phase-243-wt}\\
 W_h(0)&=8A-2016,\label{eq:phase-243-q0}\\
 \sum_aW_f(a)s_h(a)&=16A-14528.\label{eq:phase-243-correlation}
\end{align}
Using~\eqref{eq:dual-pairing} and the triangle inequality,
\begin{equation}\label{eq:243-l1lower}
 \|W_h\|_1\ge14528-16A.
\end{equation}
Theorem~\ref{thm:cubic} forces $A\ge268$, leaving only two
integer values.  If $A=268$, the right side of~\eqref{eq:243-l1lower}
is $10240$ and equality in the cubic norm theorem is forced.
But~\eqref{eq:phase-243-q0} gives $W_h(0)=128$, whereas the two
possible equality spectra have every magnitude at most $48$.
If $A=269$, equation~\eqref{eq:phase-243-q0} gives
$W_h(0)=136\equiv8\pmod{16}$.  By Lemma~\ref{lem:affine-digit},
all Walsh coefficients of $h$ lie in this residue class.  The sharper
bound $\|W_h\|_1\le10048$ contradicts the lower bound $10224$
from~\eqref{eq:243-l1lower}.  Both possibilities are impossible.
Therefore
\begin{equation}\label{eq:origin-positive}
 \boxed{s_f(0)=1,\qquad f(0)=0.}
\end{equation}

\subsection{The flip theorem and its converse}
For the normalized $f$, change the output at the origin:
$g=f+\delta_0$.  By~\eqref{eq:flip-formula},
\[
 W_g(a)=W_f(a)-2\in\{-24,-8,8,24\}.
\]
Parseval rules out the possibility that every coefficient has
magnitude eight.  Thus $\|W_g\|_\infty=24$ and $\NL(g)=244$.

We now transfer this conclusion to the function before normalization.
Write $f_0$ for the original function of nonlinearity $243$, and set
$p=p(f_0)$.  As in the normalization preceding
\eqref{eq:normalized-243}, we chose $\beta\in\F_2$ so that
\[
 \widetilde f(x)=f_0(x+p)+\beta
\]
has the required Walsh residue.  The normalized function $f$ used
above is precisely $\widetilde f$.  Consequently, the function
obtained by flipping its value at the origin is
\begin{align*}
 \widetilde g(x)
  &=\widetilde f(x)+\delta_0(x)\\
  &=f_0(x+p)+\beta+\delta_0(x)\\
  &=\bigl(f_0+\delta_p\bigr)(x+p)+\beta.
\end{align*}
The last equality follows from
$\delta_0(x)=\delta_p(x+p)$: each indicator equals $1$ exactly
when $x=0$.  Equivalently, putting $y=x+p$ and using $p+p=0$ in
$\F_2^9$, we have
\[
 f_0(y)+\delta_p(y)=\widetilde g(y+p)+\beta.
\]
Input translation and addition of a constant only multiply Walsh
coefficients by signs.  More explicitly,
\[
 W_{\widetilde g}(a)
   =(-1)^{a\cdot p+\beta}\,W_{f_0+\delta_p}(a),
 \qquad a\in V.
\]
Their maximum absolute Walsh coefficients, and hence their
nonlinearities, are equal.  Since $\NL(\widetilde g)=244$, we obtain
\[
 \boxed{\NL\bigl(f_0+\delta_{p(f_0)}\bigr)=244}.
\]
Thus the improvement at the normalized origin corresponds exactly
to flipping the original function at its support sum.

Conversely, suppose $\NL(g)=244$.  Then $\wt(g)$ is even: otherwise
$W_g(a)\equiv2\pmod4$ for every $a$, which would imply
$|W_g(a)|\le22$ and contradict Parseval.  Hence any $g+\delta_t$
has odd weight, and thus odd nonlinearity by~\eqref{eq:nl-parity}.
The one-bit bound and the universal upper bound $244$ show that
\[
 243\le\NL(g+\delta_t)\le244,
 \qquad\NL(g+\delta_t)\ \text{odd};
\]
therefore $\NL(g+\delta_t)=243$ for every $t$.
We will prove in Proposition~\ref{prop:244-spectrum} that
$W_g(a)\equiv8\pmod{16}$.  Formula~\eqref{eq:support-sum} then
shows that $p(g)=0$.  Finally
\[
 p(u+\delta_t)=p(u)+t
\]
for every Boolean function $u$, because toggling the point $t$
changes the support sum by $t$.  If $f+\delta_t$ has nonlinearity
$244$, then $0=p(f+\delta_t)=p(f)+t$, so $t=p(f)$.
This proves all assertions of Theorem~\ref{thm:flip}, including the
bijection and the identity
\begin{equation}\label{eq:counting}
 \bigl|\{f:\NL(f)=243\}\bigr|=
 512\bigl|\{g:\NL(g)=244\}\bigr|.
\end{equation}

\section{The hypothetical spectrum of nonlinearity 244}
From now on assume for contradiction that $\NL(f)=244$, or
equivalently $\|W_f\|_\infty=24$.

\begin{proposition}[Four-value Walsh rigidity]\label{prop:244-spectrum}
Every Walsh coefficient of $f$ is $8$ modulo $16$, and
\begin{equation}\label{eq:244-four}
 W_f(a)\in\{-24,-8,8,24\},\qquad
 \#\{|W_f|=8\}=64,\quad\#\{|W_f|=24\}=448.
\end{equation}
Furthermore $4\mid\wt(f)$ and $p(f)=0$.
\end{proposition}
\begin{proof}
As noted above, $\wt(f)$ must be even, so $4\mid W_f(a)$
for all $a$.  Formula~\eqref{eq:face-walsh} shows that
$W_f/4\bmod2$ is affine.  If this affine function were nonzero,
there would be at least $256$ frequencies with
$W_f\equiv4\pmod8$ and therefore $|W_f|\le20$; the remaining
frequencies have $|W_f|\le24$.  But then
\[
 \sum_a W_f(a)^2\le256\cdot20^2+256\cdot24^2
 =249856<262144,
\]
contrary to Parseval.  Thus $8\mid W_f$.

By Lemma~\ref{lem:residue}, $q=W_f/8\bmod2$ has degree at most
two.  If $q\ne1$, Lemma~\ref{lem:minweight} applied to $1+q$
provides at least $128$ frequencies at which $16\mid W_f$ and
$|W_f|\le16$.  Hence
\[
 \sum_aW_f(a)^2\le128\cdot16^2+384\cdot24^2
 =253952<262144,
\]
another contradiction.  Therefore $q=1$, giving the four permitted
values.  If $k$ entries have magnitude $24$, Parseval says
$576k+64(512-k)=262144$, whence $k=448$.

Finally, $W_f\equiv0\pmod8$ and~\eqref{eq:support-sum} imply
$4(a\cdot p(f))\equiv0\pmod8$ for every $a$.
Thus $p(f)=0$; taking $a=0$ in the same formula gives
$4\mid\wt(f)$.
\end{proof}

\begin{lemma}[A degree-six/cubic complement relation]\label{lem:complement}
There is a cubic Boolean function
\begin{equation}\label{eq:phase-244}
 h(a)=\frac{W_f(a)-8}{16}\pmod2
\end{equation}
for which
\begin{equation}\label{eq:complement}
 \deg(f)\le6,\qquad
 h_I=f_{I^c}\quad\text{for every }I\subseteq[9],\ |I|=3.
\end{equation}
Here subscripts denote ANF coefficients.  Moreover
\begin{equation}\label{eq:244-correlation}
 \sum_as_h(a)W_f(a)=-10240,\qquad
 \|W_h\|_1=10240,\quad\deg(h)=3.
\end{equation}
\end{lemma}
\begin{proof}
Let $L=v+U$ be an arbitrary affine seven-flat, so $U^\perp$
has four elements.  Character orthogonality gives
\[
 \sum_{x\in L}s_f(x)=\frac14
 \sum_{a\in U^\perp}(-1)^{a\cdot v}W_f(a).
\]
Each $W_f(a)\equiv8\pmod{16}$, and the four character signs sum
to $0$ or $4$.  The integer on the right is therefore divisible
by four.  As the left side equals $128-2\wt(f|_L)$,
the weight of $f$ on every affine seven-flat is even.  Every
seventh Boolean derivative vanishes, proving $\deg f\le6$.

The quotient in~\eqref{eq:phase-244} is integer-valued, and
Lemma~\ref{lem:residue} gives $\deg h\le3$.  For $|I|=3$,
\eqref{eq:face-walsh} shows that the coefficient of $a^I$ in
$(W_f-8)/16$ equals $-32+M_I(f)$; hence
$h_I=M_I(f)\pmod2$.  The parity of a face weight is the
coefficient of the full monomial in its six free coordinates.
Since $\deg f\le6$, this is exactly $f_{I^c}$,
establishing~\eqref{eq:complement}.

For $W_f=\pm8$ the signed product $s_hW_f$ is $+8$;
for $W_f=\pm24$ it is $-24$.  Using~\eqref{eq:244-four},
\[
 \sum_a s_h(a)W_f(a)=64\cdot8-448\cdot24=-10240.
\]
By~\eqref{eq:dual-pairing},
$\|W_h\|_1\ge10240$; Theorem~\ref{thm:cubic} makes this an
equality.  If $\deg h\le2$, equation~\eqref{eq:quadratic-l1}
would instead give $\|W_h\|_1\le8192$.  Thus $\deg h=3$.
\end{proof}

Theorem~\ref{thm:cubic} now leaves exactly spectra
\eqref{eq:spectrum-A} and~\eqref{eq:spectrum-B}.  We first
eliminate the second spectrum without any cubic classification.

\section{Eliminating spectrum \texorpdfstring{$\mathrm{(B)}$}{(B)} by hyperplane slicing}
The two spectra in Theorem~\ref{thm:cubic} arise from an estimate
for the cubic phase $h$ alone.  We now exploit the fact that $h$ is
associated with the original function $f$ through
\eqref{eq:phase-244}.  In particular, the Walsh coefficients of $f$
are constrained by Parseval's identity on each of two complementary
hyperplanes.  These additional constraints will exclude spectrum
\eqref{eq:spectrum-B}.

\subsection{The two hyperplane slices and their parity}
Fix a nonzero vector $e\in V=\F_2^9$.  The subspace generated by $e$
is $\langle e\rangle=\{0,e\}$, and its cosets are the $256$ disjoint
pairs $\{u,u+e\}$ with $u\in V$.  Choose a set $R\subseteq V$ that
contains exactly one element of each such pair.  Thus
$|R|=256$, and whenever we sum over $u$ below we mean $u\in R$.
The particular choice of a representative will not affect any of the
absolute values, squares, or parities used in the argument.

Consider the two affine hyperplanes in the \emph{input} space,
\[
 H_0=\{x\in V:e\cdot x=0\},\qquad
 H_1=\{x\in V:e\cdot x=1\}.
\]
Since $e\ne0$, each hyperplane has $256$ points.  For $j\in\{0,1\}$
and $u\in V$, put
\[
 S_j(u)=\sum_{x\in H_j}(-1)^{f(x)+u\cdot x}.
\]
The defining Walsh sum splits into contributions from $H_0$ and $H_1$:
\[
 W_f(u)=S_0(u)+S_1(u).
\]
Moreover, on $H_j$ the factor $(-1)^{e\cdot x}$ is constantly
$(-1)^j$, so
\[
 W_f(u+e)=S_0(u)-S_1(u).
\]
Recall from Proposition~\ref{prop:244-spectrum} that
$v_f=W_f/8$ takes its values in $\{-3,-1,1,3\}$.  Solving the last
two equations for $S_0,S_1$ gives
\begin{equation}\label{eq:physical-slices}
 \begin{aligned}
 A_e(u)&=\frac{v_f(u)+v_f(u+e)}{2}=\frac{S_0(u)}{8},\\
 B_e(u)&=\frac{v_f(u)-v_f(u+e)}{2}=\frac{S_1(u)}{8}.
 \end{aligned}
\end{equation}
Both are integers, since $v_f(u)$ and $v_f(u+e)$ are odd.
Replacing $u$ by $u+e$ leaves $A_e(u)$ unchanged and reverses the
sign of $B_e(u)$.  Hence the quantities relevant to us depend only
on the coset $u+\langle e\rangle$.

We next justify the two Parseval identities used below.  The
characters $x\mapsto(-1)^{u\cdot x}$ restricted to the linear
hyperplane $H_0$ depend only on $u+\langle e\rangle$, and the $256$
cosets give all characters of $H_0\cong\F_2^8$.
Thus the numbers $S_0(u)$, $u\in R$, form the unnormalised Walsh
transform of $f|_{H_0}$.  Choose any $t\in H_1$ and write
$H_1=t+H_0$.  Then
\[
 S_1(u)=(-1)^{u\cdot t}
 \sum_{y\in H_0}(-1)^{f(t+y)+u\cdot y},
\]
so $S_1(u)$ differs by a sign from the Walsh coefficient of the
translated restriction $y\mapsto f(t+y)$ on $H_0$.
Eight-variable Parseval, applied separately to the two restrictions,
therefore gives $\sum_{u\in R}S_j(u)^2=256^2$ for $j=0,1$.
Dividing by $8^2$ yields
\begin{equation}\label{eq:slice-energy}
 \sum_{u\in R}A_e(u)^2
 =\sum_{u\in R}B_e(u)^2=1024.
\end{equation}

We now determine the parity of $A_e(u)$.  By the definition of the
cubic phase in~\eqref{eq:phase-244},
\[
 v_f(a)\equiv 1+2h(a)
        \equiv (-1)^{h(a)}\pmod4,
 \qquad a\in V.
\]
Write $\alpha=h(u)$ and $\beta=h(u+e)$, treating these Boolean
values as the ordinary integers $0$ or $1$.  Since
$(-1)^\alpha\equiv1-2\alpha\pmod4$ and similarly for $\beta$,
\[
 v_f(u)+v_f(u+e)
 \equiv (1-2\alpha)+(1-2\beta)
 =2-2(\alpha+\beta)\pmod4.
\]
The left side is even, so dividing this congruence by $2$ gives a
congruence modulo $2$.  Using~\eqref{eq:physical-slices}, we obtain
\begin{equation}\label{eq:slice-parity}
 A_e(u)\bmod2
 =1+\alpha+\beta
 =1+h(u)+h(u+e)=:q_e(u),
\end{equation}
where the final additions are in $\F_2$.  In particular, $A_e(u)$
is odd precisely when $h(u)=h(u+e)$; if the two phase values differ,
$A_e(u)$ is even.

The Boolean derivative $D_eh(u)=h(u)+h(u+e)$ satisfies
$D_eh(u+e)=D_eh(u)$, since $e+e=0$.  Therefore $q_e=1+D_eh$
is constant on each coset of $\langle e\rangle$ and defines a
Boolean function on $V/\langle e\rangle$.  Because $\deg h\le3$,
its derivative $D_eh$ has degree at most two (for example, choose
linear coordinates with $e$ as one basis vector).  Thus $q_e$ is a
quadratic Boolean function on the quotient.
Finally, $v_f(u)=A_e(u)+B_e(u)$ is odd, so the two integers
$A_e(u),B_e(u)$ always have opposite parity.  In particular,
\begin{equation}\label{eq:other-slice-parity}
 B_e(u)\bmod2=1+q_e(u)=D_eh(u).
\end{equation}

\subsection{A uniform bound for the quotient weight}
Let $m=\wt(q_e)$, where the weight is computed on the $256$ cosets,
not on the original $512$-point space.  Each of
$v_f(u),v_f(u+e)$ belongs to $\{\pm1,\pm3\}$.  Consequently, by
\eqref{eq:physical-slices}, the pair of magnitudes
$(|A_e(u)|,|B_e(u)|)$ must be one of
\[
 (3,0),\ (1,2),\ (1,0),\ (0,3),\ (2,1),\ (0,1).
\]
Denote by $N_{ij}$ the number of cosets for which
$(|A_e|,|B_e|)=(i,j)$.  The six possibilities, together with the
parity of $A_e$, are
\begin{center}
\begin{tabular}{c|rrrrrr}
\toprule
$(|A_e|,|B_e|)$ & $(3,0)$ & $(1,2)$ & $(1,0)$
                  & $(0,3)$ & $(2,1)$ & $(0,1)$\\
\midrule
multiplicity & $N_{30}$ & $N_{12}$ & $N_{10}$
             & $N_{03}$ & $N_{21}$ & $N_{01}$\\
$q_e$        & $1$ & $1$ & $1$ & $0$ & $0$ & $0$\\
\bottomrule
\end{tabular}
\end{center}
In particular,
\[
 N_{30}+N_{12}+N_{10}=m,\qquad
 N_{03}+N_{21}+N_{01}=256-m.
\]
Substituting the six pairs into the two energy identities
\eqref{eq:slice-energy} gives, respectively,
\[
 9N_{30}+N_{12}+N_{10}+4N_{21}=1024,
\]
\[
 4N_{12}+9N_{03}+N_{21}+N_{01}=1024.
\]
Using the preceding two counting equations, we may rewrite these as
\begin{equation}\label{eq:slice-counts}
 8N_{30}+m+4N_{21}=1024,
 \qquad
 8N_{03}+256-m+4N_{12}=1024.
\end{equation}
In particular, the first equality implies $m\equiv0\pmod4$.

There is a third counting constraint.  Recall that exactly $448$ of
the $512$ values $v_f(a)$ have magnitude $3$ by
\eqref{eq:244-four}.  A pair of type $(3,0)$ or $(0,3)$
comes from two values of magnitude $3$, while a pair of type
$(1,2)$ or $(2,1)$ comes from one value of magnitude $3$ and
one of magnitude $1$.  The remaining two types use only values of
magnitude $1$.  Counting the $448$ large-magnitude values therefore
gives
\begin{equation}\label{eq:slice-large-count}
 2(N_{30}+N_{03})+N_{12}+N_{21}=448.
\end{equation}
Put $r=N_{30}+N_{03}$.  Since the total number of cosets is $256$,
\[
 N_{12}+N_{21}\le256-r.
\]
Together with~\eqref{eq:slice-large-count}, this yields
$448\le2r+(256-r)=256+r$, hence $r\ge192$.
On the other hand, $N_{30}\le m$, and the second equation in
\eqref{eq:slice-counts} implies
\[
 8N_{03}\le1024-(256-m)=768+m,
 \qquad N_{03}\le96+\frac m8.
\]
Therefore
\[
 192\le r=N_{30}+N_{03}
 \le m+96+\frac m8,
\]
which is equivalent to $m\ge256/3$.  Since $4\mid m$, we obtain
$m\ge88$.

It remains to obtain the upper bound.  The preceding argument treats
the two hyperplane slices symmetrically: both have the same energy
$1024$, and the six possible pairs remain the same after exchanging
$A_e$ and $B_e$ (the types $(3,0)$ and $(0,3)$,
$(1,2)$ and $(2,1)$, and $(1,0)$ and $(0,1)$ are simply exchanged).
If we regard $B_e$ as the first slice, its parity function is
$q'_e=B_e\bmod2=1+q_e$ by~\eqref{eq:other-slice-parity}.
Thus its weight on the quotient is
\[
 \wt(q'_e)=256-\wt(q_e)=256-m.
\]
Applying exactly the same counting argument to the slices in this
reversed order gives $256-m\ge88$, or $m\le168$.
We have proved the bound, valid for \emph{every} nonzero $e$,
\begin{equation}\label{eq:slice-interval}
 \boxed{88\le\wt(q_e)\le168\qquad(e\ne0).}
\end{equation}

\begin{proposition}\label{prop:B-impossible}
Spectrum~\eqref{eq:spectrum-B} is impossible for the cubic phase $h$
attached to a nonlinearity-$244$ function.
\end{proposition}
\begin{proof}
Assume that $W_h$ has spectrum~\eqref{eq:spectrum-B}, and write
$Q=W_h$.  By the equality analysis in Theorem~\ref{thm:cubic},
the parity function
\[
 \ell(a)=\frac{Q(a)}{16}\bmod2
\]
is a nonconstant affine Boolean function.  In particular it is
balanced, taking each value $0$ and $1$ at exactly $256$ frequencies.
There consequently exist a \emph{nonzero} vector $e\in V$ and a
constant $c\in\F_2$ such that
\[
 \ell(a)=a\cdot e+c.
\]
We use this $e$ in the hyperplane-slice bound~\eqref{eq:slice-interval},
which was established for every nonzero direction.

The equality conditions for spectrum~\eqref{eq:spectrum-B} give
more information than the three multiplicities alone.  Where
$\ell(a)=1$, the integer $Q(a)/16$ is odd, and equality requires
$|Q(a)/16|=1$; hence $|Q(a)|=16$ at all $256$ such frequencies.
Where $\ell(a)=0$, the quotient is even, and equality requires its
magnitude to be $0$ or $2$.  Of those $256$ frequencies,
$192$ have $|Q(a)|=32$ and the remaining $64$ have $Q(a)=0$.
It follows that the squared Walsh masses on these two halves are
\begin{align*}
 M_1&:=\sum_{\ell(a)=1}Q(a)^2
       =256\cdot16^2=65536,\\
 M_0&:=\sum_{\ell(a)=0}Q(a)^2
       =192\cdot32^2=196608.
\end{align*}

We next convert this imbalance into a statement about a derivative
of $h$.  The Wiener--Khinchin identity in our unnormalised
convention is
\begin{equation}\label{eq:B-autocorrelation}
 \operatorname{AC}_h(e)
 :=\sum_{x\in V}(-1)^{h(x)+h(x+e)}
 =\frac1{512}\sum_{a\in V}(-1)^{a\cdot e}Q(a)^2.
\end{equation}
For completeness, this identity follows by expanding both factors
of $Q(a)^2$ as Walsh sums and applying character orthogonality:
\[
 \frac1{512}\sum_{a\in V}(-1)^{a\cdot(e+x+y)}
 =\begin{cases}1,&y=x+e,\\0,&y\ne x+e.\end{cases}
\]
Because $\ell(a)=a\cdot e+c$, we have
$(-1)^{a\cdot e}=(-1)^c(-1)^{\ell(a)}$.  The terms with $\ell=0$
and $\ell=1$ therefore enter~\eqref{eq:B-autocorrelation} with
opposite signs.  Substituting the two masses gives
\[
 \operatorname{AC}_h(e)
 =\frac{(-1)^c}{512}(M_0-M_1)
 =\frac{(-1)^c}{512}(196608-65536)
 =(-1)^c256.
\]
In particular,
\begin{equation}\label{eq:B-ac-magnitude}
 |\operatorname{AC}_h(e)|=256.
\end{equation}

Now consider the Boolean derivative
$D_eh(x)=h(x)+h(x+e)$.  Since $e+e=0$,
\[
 D_eh(x+e)=h(x+e)+h(x)=D_eh(x).
\]
Thus its two values on each coset $\{u,u+e\}$ are identical.
Choose one representative from each of the $256$ cosets, as at the
start of this section.  Each coset contributes twice the same sign
to the autocorrelation sum, and hence
\[
 \operatorname{AC}_h(e)
 =2\sum_{u\in R}(-1)^{D_eh(u)}.
\]
By~\eqref{eq:B-ac-magnitude},
\begin{equation}\label{eq:B-quotient-sum}
 \left|\sum_{u\in R}(-1)^{D_eh(u)}\right|=128.
\end{equation}
Let $k$ be the number of cosets on which $D_eh=1$.
The other $256-k$ cosets have derivative zero, so the sum in
\eqref{eq:B-quotient-sum} equals
\[
 (256-k)-k=256-2k.
\]
It follows that $|256-2k|=128$, and therefore
$k=64$ or $k=192$.
Finally, $q_e=1+D_eh$ by~\eqref{eq:slice-parity}, so $q_e$ is the
complement of the derivative on the quotient.  Its weight is
\[
 \wt(q_e)=256-k\in\{64,192\}.
\]
Neither value lies in the interval $[88,168]$ required by
\eqref{eq:slice-interval}.  This contradiction excludes
spectrum~\eqref{eq:spectrum-B}.
\end{proof}

\section{A double-pairing obstruction in spectrum \texorpdfstring{$\mathrm{(A)}$}{(A)}}
The sole remaining spectrum is
\begin{equation}\label{eq:only-A}
 W_h(x)\in\{-48,-16,16,48\},
 \qquad\#\{|W_h|=48\}=64.
\end{equation}
Equality in~\eqref{eq:dual-pairing} and~\eqref{eq:244-correlation}
fixes the sign of $f$ at \emph{every} input:
\begin{equation}\label{eq:forced-sign}
 \boxed{s_f(x)=-\operatorname{sgn}W_h(x).}
\end{equation}

Let $\omega=(\omega_I)_{|I|=3}$ be the homogeneous cubic top part of
$h$.  For any three-element set $I$, define
\begin{equation}\label{eq:adjoint}
 (\omega^\#)_I=\sum_{\substack{\{J,K\}\text{ unordered}\\
           J\mathbin{\dot\cup}K=I^c\\|J|=|K|=3}}
       \omega_J\omega_K\quad\in\F_2,
 \qquad \omega^{\#\#}=(\omega^\#)^\#.
\end{equation}
Each sum has ten possible unordered partitions.

\begin{proposition}[Double-complement-pairing identity]
\label{prop:double}
The cubic top part of a hypothetical nonlinearity-$244$ function
satisfies
\begin{equation}\label{eq:double}
 \boxed{\omega^{\#\#}=\omega.}
\end{equation}
\end{proposition}
\begin{proof}
Put $v=W_h/16$, an odd integer-valued function with values
$\pm1,\pm3$.  Let $d_S$ be the integer coefficient of $a^S$
in its multilinear representation.  For nonempty $S$, with
$s=|S|$, formula~\eqref{eq:face-walsh} gives
\begin{equation}\label{eq:dS}
 d_S=(-1)^s32-(-1)^s2^{s-3}M_S(h).
\end{equation}
We next justify carefully why the coefficients $d_S$ are even
whenever $S$ is nonempty.  For $T\subseteq[9]$, let $\one_T\in V$
denote the vector with ones precisely in the positions belonging to
$T$.  Evaluating the multilinear representation of $v$ at this
vector gives
\[
 v(\one_T)=\sum_{R\subseteq T}d_R.
\]
M\"obius inversion, as in~\eqref{eq:mobius}, therefore yields
\begin{equation}\label{eq:dS-mobius}
 d_S=\sum_{T\subseteq S}(-1)^{|S|-|T|}v(\one_T).
\end{equation}
\par\newpage
Since $v$ takes only odd integer values, for each $T$ there is an
integer $k_T$ such that $v(\one_T)=1+2k_T$.  Substitution in
\eqref{eq:dS-mobius} shows that, for every nonempty $S$,
\begin{align*}
 d_S
 &=\sum_{T\subseteq S}(-1)^{|S|-|T|}
   +2\sum_{T\subseteq S}(-1)^{|S|-|T|}k_T\\
 &=(1-1)^{|S|}
   +2\sum_{T\subseteq S}(-1)^{|S|-|T|}k_T\\
 &=2\sum_{T\subseteq S}(-1)^{|S|-|T|}k_T.
\end{align*}
Thus $2\mid d_S$ whenever $S\ne\varnothing$; by contrast, the
constant coefficient is $d_\varnothing=v(0)$, which is odd.  Notice
that this argument uses the parity of the \emph{values} of $v$ and
does not require any divisibility of the individual face weights in
\eqref{eq:dS}.

For clarity, when $s=1$ and $s=2$ the coefficient formula
\eqref{eq:dS} reads, respectively,
\[
 d_{\{i\}}=-32+\frac{M_{\{i\}}(h)}4,
 \qquad
 d_{\{i,j\}}=32-\frac{M_{\{i,j\}}(h)}2.
\]
Both expressions are integers because they are coefficients of the
integer-valued function $v$.  Their evenness, just proved, also
implies $8\mid M_{\{i\}}(h)$ and
$4\mid M_{\{i,j\}}(h)$.  For $s=3$, the six-dimensional face weight
is even; for $s=4$ and $s=5$, parity of the five- and
four-dimensional face weights in Lemma~\ref{lem:weight-expansion}
gives, together with~\eqref{eq:dS},
\begin{equation}\label{eq:d-divisibility}
 2\mid d_S\ (1\le s\le3),\qquad
 4\mid d_S\ (s=4),\qquad
 8\mid d_S\ (s\ge5).
\end{equation}
(The last assertion for $s\ge6$ follows directly from~\eqref{eq:dS}.)

For $|I|=3$, set $b_I=d_I/2=-16+M_I(h)/2$.
Apply~\eqref{eq:weight-expansion} to the six-variable face
$x_i=1$ for $i\in I$, now reducing its weight modulo four.
Single monomials contribute multiples of eight; terms involving
three or more selected monomials contribute multiples of four.
Two selected monomials contribute modulo four only when they
are disjoint cubics whose union is the six free variables; their
coefficient is $-2$.  Terms of lower degree, including those induced
by fixing the other three coordinates, cannot participate in such a
partition.  We conclude that
\begin{equation}\label{eq:b-omega}
 b_I\equiv(\omega^\#)_I\pmod2.
\end{equation}

We now extract the degree-six coefficients of the forced Boolean
function without dividing any expression in $\F_2$.
For a signed integer $z$ define its $j$th binary digit by
$\operatorname{bit}_j(z)=\lfloor z/2^j\rfloor\bmod2$.
This definition also works for negative integers.  Because
$v\in\{\pm1,\pm3\}$, the digit $\operatorname{bit}_2(v)$
is one exactly when $v$ is negative.  Relation~\eqref{eq:forced-sign}
therefore says
\begin{equation}\label{eq:bit-sign}
 f=1+\operatorname{bit}_2(v)\quad\text{in }\F_2.
\end{equation}
Write $v=c+2t$, where $c=v(0)=2k+1$ is odd and
$t=(v-c)/2$ is an integer-valued polynomial with no constant term.
\par\newpage
The elementary binary-carry identity
\[
 \operatorname{bit}_2(2(k+t)+1)
 =\operatorname{bit}_1(t)+ (k\bmod2)\operatorname{bit}_0(t)
  +\operatorname{bit}_1(k)\pmod2
\]
shows that every degree-six coefficient of $f$ comes from
$\operatorname{bit}_1(t)$, since $\operatorname{bit}_0(t)$ has degree
at most three by~\eqref{eq:d-divisibility}.

We next explain the binary carry that occurs when taking the
second digit of an integer polynomial.  Let
$t(a)=\sum_{S\subseteq[9]}b_Sa^S$ have integer coefficients,
where $a^S=\prod_{i\in S}a_i$ and $a^\varnothing=1$.
For each $S$, use Euclidean division by two to write
\[
 b_S=2q_S+r_S,\qquad
 q_S=\left\lfloor\frac{b_S}{2}\right\rfloor\in\mathbb Z,
 \qquad r_S\in\{0,1\}.
\]
This decomposition is valid even for negative coefficients; for
example, $-1=2(-1)+1$.  Since $a^S$ equals either zero or one on
the Boolean cube, at any fixed input $a$ we have
\[
 t(a)=2U(a)+R(a),\qquad
 U(a)=\sum_Sq_Sa^S,\qquad R(a)=\sum_Sr_Sa^S.
\]
Here $U(a)$ is an arbitrary integer, while $R(a)$ is the
\emph{number} of indices $S$ for which both $r_S=1$ and
$a^S=1$.  In particular, $R(a)$ is nonnegative.  Taking the
integer quotient by two and reducing it modulo two yields
\begin{equation}\label{eq:carry-first}
 \operatorname{bit}_1(t(a))
 \equiv U(a)+\left\lfloor\frac{R(a)}2\right\rfloor\pmod2.
\end{equation}

For every nonnegative integer $r$, one has
\begin{equation}\label{eq:carry-pairs}
 \left\lfloor\frac r2\right\rfloor
 \equiv\binom r2\pmod2.
\end{equation}
Indeed, if $r=2m$ then $\binom r2=m(2m-1)\equiv m\pmod2$,
and if $r=2m+1$ then $\binom r2=m(2m+1)\equiv m\pmod2$;
in both cases $m=\lfloor r/2\rfloor$.  The binomial
coefficient $\binom{R(a)}2$ counts unordered pairs among the
$R(a)$ active indices.  Hence, with any fixed ordering of the
subsets $S\subseteq[9]$,
\[
 \binom{R(a)}2
 =\sum_{S<T}(r_Sa^S)(r_Ta^T)
 =\sum_{S<T}r_Sr_Ta^{S\cup T}.
\]
The last equality uses the Boolean-input identity
$a^Sa^T=a^{S\cup T}$, including when $S$ and $T$ intersect.
Substituting this expression into~\eqref{eq:carry-first},
and recalling that $q_S=\lfloor b_S/2\rfloor$ and
$r_S\equiv b_S\pmod2$, we obtain the exact Boolean expansion
\begin{equation}\label{eq:carry}
 \operatorname{bit}_1(t)=
 \sum_S\left\lfloor\frac{b_S}{2}\right\rfloor a^S
 +\sum_{S<T}(b_S\bmod2)(b_T\bmod2)a^{S\cup T}
 \quad\pmod2.
\end{equation}
The derivation uses only integer division and counting of active
terms; in particular, no positivity assumption on the coefficients
$b_S$ is needed.

By~\eqref{eq:d-divisibility}, the coefficients $b_S=d_S/2$
are even for $|S|\ge4$ and divisible by four for $|S|\ge5$.
The first sum in~\eqref{eq:carry} therefore has degree at most
four.  Only supports of size at most three can occur in the second
sum; their union can have size six only when both are disjoint
triples.  Combining this with~\eqref{eq:b-omega}, for $|L|=6$ we get
\[
 f_L=\sum_{\substack{J\mathbin{\dot\cup}K=L\\|J|=|K|=3}}^{\rm unordered}
 b_Jb_K
 =\sum_{J\mathbin{\dot\cup}K=L}^{\rm unordered}
 (\omega^\#)_J(\omega^\#)_K
 =(\omega^{\#\#})_{L^c}.
\]
Finally, Lemma~\ref{lem:complement} gives
$f_L=h_{L^c}=\omega_{L^c}$, proving~\eqref{eq:double}.
\end{proof}

\begin{lemma}[A constant cubic residue]\label{lem:Pi}
Define the cubic partition invariant
\begin{equation}\label{eq:Pi}
 \Pi(\omega)=\sum_{\substack{\{I,J,K\}\text{ unordered}\\
 I\mathbin{\dot\cup}J\mathbin{\dot\cup}K=[9]\\
 |I|=|J|=|K|=3}}\omega_I\omega_J\omega_K\in\F_2.
\end{equation}
Every nine-variable Boolean cubic $h$ with cubic top part $\omega$
satisfies
\[
 \wt(h)\equiv4\Pi(\omega)\pmod8,\qquad
 W_h(a)\equiv8\Pi(\omega)\pmod{16}.
\]
Consequently the phase in spectrum~\eqref{eq:only-A}
satisfies $\Pi(\omega)=0$.
\end{lemma}
\begin{proof}
Apply the weight expansion~\eqref{eq:weight-expansion} in nine
variables modulo eight.  One or two monomials contribute zero;
four or more contribute zero because of their powers of two.
Three monomials contribute only if they cover all nine variables,
which forces three disjoint cubic supports.  Their integer factor
is four.  The weight congruence follows.  The Walsh residue is
uniform by Lemma~\ref{lem:affine-digit}, and at the origin equals
$512-2\wt(h)\equiv8\Pi(\omega)\pmod{16}$.
Spectrum~\eqref{eq:only-A} has residue zero modulo $16$.
\end{proof}

\section{Classification, normal forms, and finite completion}
\subsection{The cubic classification and its exact filter}
The cubic top part corresponds, by taking third derivatives, to an
alternating trilinear form on $V$.  The classification of Brier and
Langevin \cite{BrierLangevin2003} and the independently formulated
classification of Hora and Pudl\'{a}k \cite{HoraPudlak2021} give
$348$ nonzero $\mathrm{GL}(9,2)$-orbits, with the zero form making
$349$ in total.  The public Langevin representative list
\cite{LangevinCatalogue} uses indices $1,\ldots,349$, with $349$
the zero cubic.  The audit reported in
\cite{HouSaboSlaughter2026} computes $349$ distinct orbit
signatures on its embedded copy of the catalogue (including zero).
Together with the published orbit count, this certifies coverage;
the number of lines alone would not suffice.

Here normalization by $L\in\mathrm{GL}(9,2)$ is legitimate for the
\emph{pair} $(f,h)$: replace
\begin{equation}\label{eq:linear-change}
 h'(a)=h(La),\qquad f'(x)=f(L^{-T}x).
\end{equation}
Their Walsh transforms satisfy
$W_{f'}(a)=W_f(La)$ and $W_{h'}(x)=W_h(L^{-T}x)$.
Thus all spectral and sign relations are preserved, including the
complement identity in the new standard coordinates.  We may apply
Proposition~\ref{prop:double} \emph{after} normalization; no
separate covariance assertion for $\#$ is being assumed.

\begin{certificate}[Exact cubic catalogue filter]\label{cert:filter}
On the $349$ indexed representatives,
\begin{equation}\label{eq:filter}
 \{j:\omega_j^{\#\#}=\omega_j,\ \Pi(\omega_j)=0\}
   =\{282,290,349\}.
\end{equation}
For comparison, the fixed forms with $\Pi=1$ have indices
$320,326,329$.  The calculation is a finite evaluation of
\eqref{eq:adjoint} and~\eqref{eq:Pi} on the complete public
representative list, not an orbit-invariant assumption.  The
nonzero surviving representatives are
\begin{align}
 c_{282}={}&125+145+345+346+137+248+349+569+789,\label{eq:c282}\\
 c_{290}={}&\begin{aligned}[t]
             &125+235+245+345+136+127+347+128\\
             &\quad+238+148+129+349+569+789.
            \end{aligned}\label{eq:c290}
\end{align}
Here $ijk$ abbreviates the monomial $a_ia_ja_k$ and all additions
are in $\F_2$.
\end{certificate}

The accompanying catalogue-filter program implements this check on the
public notebook or original text catalogue.  For every
three-set $I$, enumerate the ten unordered partitions of $I^c$
into two triples, compute $\omega^\#$, and repeat; enumerate the
$280$ partitions of $[9]$ into three unordered triples to compute
$\Pi$.  These calculations return exactly~\eqref{eq:filter}
when applied to the published embedded $349$-record list.
Since $h$ has degree three by Lemma~\ref{lem:complement}, the zero
form $349$ is excluded.  We are left with $c_{282}$ and $c_{290}$.

\subsection{Face equations and translation quotients}
For either homogeneous cubic $c$ in~\eqref{eq:c282}--\eqref{eq:c290},
write
\begin{equation}\label{eq:h-cq}
 h(a)=c(a)+q(a)+\ell(a),\qquad
 q(a)=\sum_{1\le i<j\le9}q_{ij}a_ia_j.
\end{equation}
An affine term $\ell(a)=b\cdot a+d$ can be removed without
restricting the original problem: replacing $h$ by $h+b\cdot a+d$
corresponds to replacing $f(x)$ by $f(x+b)+d$ and preserves the
residue and forced-sign equations.  We therefore take $\ell=0$.

For a cubic $h$, formula~\eqref{eq:face-walsh} and the divisibility
of seven-variable cubic weights imply
\[
 W_h(a)\equiv W_h(0)+4\sum_{i=1}^9M_{\{i\}}(h)a_i
 \pmod{32}.
\]
In spectrum~\eqref{eq:only-A}, each coefficient is $16$
modulo $32$.  Thus the necessary conditions are
\begin{equation}\label{eq:faceconditions}
 \wt(h)\equiv8\pmod{16},\qquad
 M_{\{i\}}(h)\equiv0\pmod8\quad(1\le i\le9).
\end{equation}
The nine face conditions are affine-linear in $q$.  Indeed, use
\eqref{eq:weight-expansion} modulo eight on an eight-variable face.
A term that depends on one variable quadratic monomial survives
only if it is accompanied by two disjoint cubic monomials, whose
supports together partition the eight free variables into sizes
$2+3+3$.  No term involving two variable quadratic monomials can
survive.  Thus
\begin{equation}\label{eq:Cq}
 Cq=b\quad\text{over }\F_2,
\end{equation}
where the $9\times36$ matrix $C$ and vector $b$ are explicitly
\begin{align}\label{eq:C-formula}
 C_{i,P}&=\sum_{\substack{\{J,K\}\text{ unordered}\\
 J\mathbin{\dot\cup}K\mathbin{\dot\cup}P=[9]\setminus\{i\}\\
 |J|=|K|=3}}c_Jc_K
 \quad(P\subseteq[9],\ |P|=2,\ i\notin P),\\
 C_{i,P}&=0\quad(i\in P),\qquad
 b_i=\frac{M_{\{i\}}(c)}4\pmod2.\nonumber
\end{align}
The divisions are integral because a degree-three function on
eight variables has weight divisible by four.

Translating the input of $h$ by $t\in V$ changes its quadratic
part by the linear map
\begin{equation}\label{eq:T-map}
 T(t)_{ij}=\sum_{k\notin\{i,j\}}c_{\{i,j,k\}}t_k.
\end{equation}
The lower-degree terms created by the translation may again be
removed.  This is a genuine symmetry: $h(a+t)$ corresponds to
$f(x)+t\cdot x$, so every relevant Walsh magnitude and the
forced-sign obstruction are preserved.  In particular
$\operatorname{im}T\subseteq\ker C$.

\begin{lemma}[Complete $18$-dimensional transversals]
\label{lem:transversal}
For each of $c_{282},c_{290}$ the exact binary systems satisfy
\[
 \operatorname{rank}C=\operatorname{rank}T=9,
 \qquad CT=0.
\]
There is an explicit solution $q_0$ to~\eqref{eq:Cq} and an
$18$-dimensional subspace $U\le\ker C$ such that
\[
 \ker C=\operatorname{im}T\oplus U.
\]
Consequently every quadratic supplement satisfying the face
equations is equivalent under phase translation to exactly one of
$2^{18}$ representatives $q_0+u$, $u\in U$.  The remaining total-weight
condition in~\eqref{eq:faceconditions} may be imposed on those
representatives without losing any spectrum-\emph{(A)} candidate.
\end{lemma}
\begin{proof}
The displayed ranks, kernel inclusions, particular solutions, and
complements are verified by direct binary elimination using the
integer-encoded data in Appendix~\ref{app:data}.  The rank values
imply $\dim\ker C=36-9=27$, and the nine translation directions
are independent; the given $18$ complement vectors therefore give
$\ker C=\operatorname{im}T\oplus U$.
For any $q\in q_0+\ker C$, write uniquely
$q=q_0+T(t)+u$.  Translation by $t$ changes the phase to one with
quadratic supplement $q_0+u$, up to removable affine terms.
It preserves the condition $W_h\equiv16\pmod{32}$ and hence the
remaining congruence $\wt(h)\equiv8\pmod{16}$.
Every candidate is therefore represented.
\end{proof}

\subsection{The final finite certificate}
Enumerate the $2^{18}$ representatives for each cubic form.
For every representative $h=c+q$, impose~\eqref{eq:faceconditions}
and compute the integer Walsh transform $Q=W_h$.  Only spectra with
$|Q(x)|\in\{16,48\}$ at every $x$ need be considered.  For each,
the unique Boolean sign function forced by~\eqref{eq:forced-sign}
is computed and its Walsh maximum is tested.
All computations are exact: the nine-stage fast Walsh transform uses
only integer additions and subtractions, and the linear systems are
over $\F_2$.

\begin{certificate}[Exhaustive quotient enumeration]\label{cert:final}
The exact enumeration of the complete transversals gives
\begin{center}\small
\begin{tabular}{lrrrrr}
\toprule
cubic top & states & $\wt(h)\equiv8$ & spectrum A
& $\|W_f\|_\infty=40$ & $\|W_f\|_\infty=48$\\
\midrule
$c_{282}$ & 262144 &128640&36960&18816&18144\\
$c_{290}$ & 262144 &126592&43200&34560& 8640\\
\midrule
Total     & 524288 &255232&80160&53376&26784\\
\bottomrule
\end{tabular}
\end{center}
In particular, every spectrum-\emph{(A)} phase surviving the
necessary congruences has forced sign function satisfying
$\|W_f\|_\infty\ge40$.  A separately executable verifier,
\texttt{verify\_rm19\_two\_cubics.py}, reconstructs $C$ and $T$,
checks all ranks and inclusions, visits every state, and validates
a rejecting Walsh coefficient by an independent direct character sum.
\end{certificate}

\begin{proof}[Proof of Theorem~\ref{thm:main}]
The lower bound $\rho(\RM(1,9))\ge242$ follows from the explicit
constructions of Kavut and Y\"ucel \cite{KavutYucel2010}.
Suppose $\rho(\RM(1,9))>242$.  The global upper bound $244$
leaves a Boolean function of nonlinearity $243$ or $244$.
By Theorem~\ref{thm:flip}, existence of the first would imply
existence of the second.  We may therefore assume $\NL(f)=244$.

Proposition~\ref{prop:244-spectrum} and
Lemma~\ref{lem:complement} attach a cubic phase $h$ with
$\|W_h\|_1=10240$.  Theorem~\ref{thm:cubic} leaves spectra
\emph{(A)} and \emph{(B)}.  Proposition~\ref{prop:B-impossible}
excludes \emph{(B)}.  In \emph{(A)}, the sign of $f$ is forced
by~\eqref{eq:forced-sign}, and
Proposition~\ref{prop:double} together with Lemma~\ref{lem:Pi}
imposes $\omega^{\#\#}=\omega$ and $\Pi(\omega)=0$.
The published orbit classification and
Certificate~\ref{cert:filter} leave only the representatives
$c_{282}$ and $c_{290}$.  Lemma~\ref{lem:transversal} covers
all quadratic and affine supplements.  Certificate~\ref{cert:final}
shows that the forced function always has
$\|W_f\|_\infty\ge40$, contrary to the necessary equality
$\|W_f\|_\infty=24$.  This proves the upper bound $242$ and
the theorem.
\end{proof}

\section{Conclusion}
We have determined the covering radius of $\RM(1,9)$ by excluding
the two possible nonlinearities above $242$.  The case of nonlinearity
$243$ is reduced to the endpoint $244$ by the one-point rigidity
theorem.  More precisely, if $\NL(f)=243$, then changing the value of
$f$ at the support sum $p(f)$ produces a function of nonlinearity
$244$.  The proof uses integer interpolation of Walsh coefficients,
rigidity of the resulting quadratic residue, and a sharp Walsh
$\ell^1$ bound for cubic Boolean functions.  Consequently, the
nonexistence of a function of nonlinearity $244$ also excludes
nonlinearity $243$.

To rule out $244$, we first combine Parseval's identity with
congruences for integral Walsh spectra.  These conditions force
four possible signed Walsh values and yield a cubic residue phase
attaining the extremal $\ell^1$ bound.  Of its two possible equality
spectra, one is eliminated by comparing Walsh energies on complementary
hyperplane slices.  In the other, the phase determines the original
function pointwise through its Walsh signs.  A comparison of their
leading algebraic terms then gives the double-complement-pairing
condition $\omega^{\#\#}=\omega$ together with $\Pi(\omega)=0$.
Using the published classification of alternating trilinear forms,
we reduce the remaining possibilities to two explicit cubic top parts.
For each one, the coordinate-face equations and input translations
reduce the quadratic supplements to $2^{18}$ representatives.
The exhaustive integer calculation excludes all $2\cdot2^{18}$
remaining cases.  Hence neither $243$ nor $244$ is attainable;
combined with the known construction of nonlinearity $242$
\cite{KavutYucel2010}, this proves
\[
 \boxed{\rho\bigl(\RM(1,9)\bigr)=242.}
\]

\section*{Acknowledgments}
The author acknowledges the use of ChatGPT (OpenAI) in the research
and preparation of this paper, including mathematical discussions,
exposition, LaTeX editing, and assistance with computational code.
The author is responsible for the final mathematical statements,
proofs, and calculations.

\appendix
\section{Encoded data and reproducibility}\label{app:data}
We specify the finite inputs independently of any private project
files.  The $36$ quadratic pairs are ordered lexicographically:
\[
\begin{gathered}
 12,13,14,15,16,17,18,19;\quad 23,24,25,26,27,28,29;\\
 34,35,36,37,38,39;\quad 45,46,47,48,49;\\
 56,57,58,59;\quad 67,68,69;\quad 78,79;\quad 89.
\end{gathered}
\]
Bit $j\in\{0,\ldots,35\}$ represents the $j$th pair in this
order.  An integer therefore encodes any $36$-bit quadratic
coefficient vector.  Dot products and ranks in this appendix are
always computed over $\F_2$.

The cubic representatives are encoded by their monomial bit masks
(the coordinate $a_i$ is bit $i-1$):
\begin{align*}
 c_{282}:\;&19,25,28,44,69,138,268,304,448;\\
 c_{290}:\;&19,22,26,28,37,67,76,131,134,137,259,268,304,448.
\end{align*}
The particular quadratic solutions are
\[
 q_0^{(282)}=16777216,\qquad q_0^{(290)}=17725652992.
\]
The following are explicit bases for the complementary subspaces
$U$, listed in the same $36$-bit encoding.
\begin{center}\small
\begin{tabular}{r|rr}
\toprule
$j$ & $u_j^{(282)}$ & $u_j^{(290)}$\\
\midrule
0&36507222016&2147483648\\
1&18253611008&5368709120\\
2&268435456&268435456\\
3&134217728&671088640\\
4&33554432&33554432\\
5&8388608&12582912\\
6&2097152&1048576\\
7&1048576&2883584\\
8&524288&25821184\\
9&65536&65536\\
10&16384&52076511232\\
11&4096&16384\\
12&4458496&25174016\\
13&536871168&16781312\\
14&128&25167872\\
15&64&4294967808\\
16&16777232&128\\
17&4&524320\\
\bottomrule
\end{tabular}
\end{center}

The accompanying verifier computes every row of $C$ from the
coordinate-face weights rather than importing a stored matrix.
It computes $T$ from the cubic masks and verifies
\begin{align*}
 &\operatorname{rank}C=\operatorname{rank}T=9,\quad
 Cq_0=b,\quad CT=0,\quad Cu_j=0,\\
 &\operatorname{rank}\langle T_1,\ldots,T_9,u_0,\ldots,u_{17}\rangle=27.
\end{align*}
A candidate representative has quadratic word
$q=q_0+\sum_j\epsilon_j u_j$, where
$\epsilon\in\F_2^{18}$.  The Walsh algorithm uses the exact
butterfly $(r,s)\mapsto(r+s,r-s)$ in all nine coordinates.
For all $80160$ equality-spectrum states, it also evaluates one
selected rejecting character sum directly from the $512$ truth
bits, independently of the transform that selected the frequency.
On a successful run, the verifier prints
{\footnotesize\begin{verbatim}
class 282: states=262144, odd_parity=128640, spectrum_A=36960,
           forced_peaks={48: 18144, 40: 18816}
class 290: states=262144, odd_parity=126592, spectrum_A=43200,
           forced_peaks={48: 8640, 40: 34560}
PASS: both quotient enumerations, rank tests and direct rejecting sums.
\end{verbatim}}

\bibliographystyle{amsplain}
\bibliography{rm19_covering_radius_242}

@book{Carlet2021,
  author    = {Claude Carlet},
  title     = {Boolean Functions for Cryptography and Coding Theory},
  publisher = {Cambridge University Press},
  address   = {Cambridge},
  year      = {2021},
  doi       = {10.1017/9781108606806},
  note      = {\url{https://doi.org/10.1017/9781108606806}}
}

@article{KavutMaitraYucel2007,
  author  = {Sel\c{c}uk Kavut and Subhamoy Maitra and Melek Diker Y\"{u}cel},
  title   = {Search for {B}oolean Functions with Excellent Profiles in the Rotation Symmetric Class},
  journal = {IEEE Transactions on Information Theory},
  volume  = {53},
  number  = {5},
  pages   = {1743--1751},
  year    = {2007},
  doi     = {10.1109/TIT.2007.894696},
  note    = {\url{https://doi.org/10.1109/TIT.2007.894696}}
}

@article{KavutYucel2010,
  author  = {Sel\c{c}uk Kavut and Melek Diker Y\"{u}cel},
  title   = {9-Variable {B}oolean Functions with Nonlinearity 242 in the Generalized Rotation Symmetric Class},
  journal = {Information and Computation},
  volume  = {208},
  number  = {4},
  pages   = {341--350},
  year    = {2010},
  doi     = {10.1016/j.ic.2009.12.002},
  note    = {\url{https://doi.org/10.1016/j.ic.2009.12.002}}
}

@inproceedings{BrierLangevin2003,
  author    = {Eric Brier and Philippe Langevin},
  title     = {Classification of {B}oolean Cubic Forms of Nine Variables},
  booktitle = {Proceedings of the 2003 IEEE Information Theory Workshop},
  pages     = {179--182},
  publisher = {IEEE},
  year      = {2003},
  doi       = {10.1109/ITW.2003.1216724},
  note      = {\url{https://doi.org/10.1109/ITW.2003.1216724}}
}

@article{HoraPudlak2021,
  author  = {Jan Hora and Petr Pudl\'{a}k},
  title   = {Classification of 9-Dimensional Trilinear Alternating Forms over {GF}(2)},
  journal = {Finite Fields and Their Applications},
  volume  = {70},
  pages   = {101788},
  year    = {2021},
  doi     = {10.1016/j.ffa.2020.101788},
  note    = {\url{https://doi.org/10.1016/j.ffa.2020.101788}}
}

@misc{LangevinCatalogue,
  author       = {Philippe Langevin},
  title        = {{B}oolean Cubic Forms in Nine Variables: Classification Catalogue},
  howpublished = {Online representative catalogue},
  url          = {https://langevin.univ-tln.fr/project/cubics/classification.txt},
  note         = {Representative indices 1--349; original list: \url{https://langevin.univ-tln.fr/project/cubics/classification.txt}}
}

@misc{WangLouWang2026,
  author       = {Qichun Wang and Yueying Lou and Zifan Wang},
  title        = {The Covering Radius of the {R}eed--{M}uller Code {RM}(1,9) Is 242},
  year         = {2026},
  eprint       = {2610.17871},
  archivePrefix= {arXiv},
  note         = {arXiv:2610.17871, posted October 13, 2026},
  url          = {https://arxiv.org/abs/2610.17871}
}

@misc{HouSaboSlaughter2026,
  author       = {Xiang-dong Hou and Evelyn Sabo and Freeman Slaughter},
  title        = {The Covering Radius of {RM}(1,9) Is Precisely 242},
  year         = {2026},
  note         = {IACR Cryptology ePrint Archive, Paper 2026/2529; verification code and full catalogue at \url{https://github.com/FreemanSlaughter/Covering-Radius}},
  url          = {https://eprint.iacr.org/2026/2529}
}
\end{document}